\documentclass[12pt, reqno]{amsart}
\usepackage{amsfonts}
\usepackage{bbm}
\usepackage{amscd,amsfonts}
\usepackage{amssymb, eucal, amsfonts, amsmath, xypic, latexsym}
\usepackage{pifont}
\usepackage{mathrsfs,color}
\usepackage{amsthm,indentfirst,bm,fancyhdr,dsfont}
\usepackage{graphicx}
\usepackage[all]{xy}
\usepackage[CJKbookmarks=true]{hyperref}

\usepackage{mathrsfs}
\usepackage{amsmath}
\usepackage{amssymb}
\usepackage{hyperref}

\newtheorem{theorem}{Theorem}[section]
\newtheorem{lemma}[theorem]{Lemma}
\newtheorem{prop}[theorem]{Proposition}
\newtheorem{proposition}[theorem]{Proposition}
\newtheorem{corollary}[theorem]{Corollary}
\theoremstyle{definition}

\newtheorem{defn}[theorem]{Definition}
\newtheorem{definition}[theorem]{Definition}

\newtheorem{observation}[theorem]{Observation}

\newtheorem{remark}[theorem]{Remark}

\newtheorem{question}[theorem]{Question}

\newcommand{\blue}[1]{{\color{blue}#1}}
\newcommand{\red}[1]{{\color{red}#1}}
\newcommand{\green}[1]{{\color{green}#1}}

\numberwithin{equation}{section}

\def\ggg{\mathfrak{g}}

\def\co{\mathcal{O}}

\def\ggg{\mathfrak{g}}

\def\ttt{\mathfrak{t}}

\def\bba{\mathbb{A}}

\def\bbz{\mathbb{Z}}

\def\bk{\mathbf{k}}

\def\bbr{\mathbb{R}}
\def\bbG{\mathbb{G}}

\def\bG{{\mathbf{G}}}

\def\bo{{\bar 1}}
\def\bz{{\bar 0}}
\def\ev{{\text{ev}}}

\def\Lie{\text{Lie}}
\def\ad{\textsf{ad}}
\def\ker{\text{Ker}}

\def\id{\mathsf{id}}

\def\GRS{\text{SGRS }}

\newcommand{\dha}{{H}\kern -.8em\mathcal{H}}

\def\salg{\mathfrak{salg}_\bk}

\def\sfr{\textsf{r}}

\def\bbU{\mathbb{U}}

\def\scrw{\mathscr{W}}

\def\sfp{\textsf{p}}

\def\frakr{\mathfrak{R}}

{\vskip-\lastskip\medskip
  \noindent
  {\em #1.}\enspace
 }
{\qed\par\medskip
}

\newcommand{\cO}{\mathcal{O}}
\newcommand{\lra}{\longrightarrow}
\newcommand{\beq}{\begin{equation}}
\newcommand{\eeq}{\end{equation}}
\newcommand{\fg}{\mathfrak{g}}

\begin{document}

\title[Basic quasi-reductive supergroups]
{Basic quasi-reductive root data and supergroups}
\author{Rita Fioresi and Bin Shu}
\address{University of Bologna
via San Donato 15,
40127 Bologna, Italy}\email{rita.fioresi@unibo.it}

\address{School of Mathematical Sciences, Ministry of Education Key Laboratory of Mathematics and Engineering Applications \& Shanghai Key Laboratory of PMMP,  East China Normal University, NO. 500 Dongchuan Rd., Shanghai 200241, China} \email{bshu@math.ecnu.edu.cn}
%\address{School of Mathematical Sciences, East China Normal University,  500 %Dongchuan Rd., Shanghai 200241,  China} \email{bshu@math.ecnu.edu.cn}


 \thanks{The research of R.F. was supported by CaLISTA CA 21109,
{MSCA-SE CaLIGOLA, Project ID: 101086123, MSCA-DN CaLiForNIA, Project ID: 101119552,
GNSAGA-Indam, PNRR MNESYS, PNRR National Center for
HPC, Big Data and Quantum Computing CUP J33C22001170001, PNNR
SIMQuSEC CUP J13C22000680006.}
The research of B.S is supported by the National Natural Science Foundation of China
(Grant No. 12071136)
% and 12271345),
and by Science and Technology Commission of Shanghai Municipality (No. 22DZ2229014). }

\begin{abstract}
  {We investigate  pairs $(G,Y)$, where $G$ is a reductive algebraic group and
  $Y$ a purely-odd $G$-superscheme, asking when a pair corresponds to
a quasi-reductive algebraic supergroup $\bbG$, that is,
$\bbG_\ev$ is isomorphic to $G$, and the quotient $\bbG\slash \bbG_\ev$
is $G$-equivariantly  isomorphic to $Y$.
We prove that, if $Y$ satisfies certain conditions (basic quasi-reductive root data),
  then the question has a positive answer given by an existence and uniqueness theorem.
  {The corresponding supergroups are said to be  basic quasi-reductive, which can be classified, up to isogeny.}
We then  decide the structure of  connected quasi-reductive algebraic supergroups
provided that:
(i) the root system does not contain $0$;
(ii) $\ggg:=\Lie(\bbG)$ admits a nondegenerate even super-symmetric bilinear form.
{(iii) all odd reflections are invertible.
Remarkably, those supergroups are exactly basic quasi-reductive supergroups of monodromy type. }

}
\end{abstract}

\maketitle

\section{Introduction}

A quasi-reductive supergroup $\bbG$ is a supergroup whose
underlying ordinary group $\bbG_\ev$ is reductive. They have been the object of study
of many papers including, for the past \cite{Ser} and refs therein, while
more recently \cite{SS}, \cite{MZ} and refs. therein. Though they represent
a natural generalization of basic supergroups, that is supergroups
with basic Lie superalgebras \cite{Kac}, \cite{FG}, %which are well known (see \cite{FG}),
they appear to be more elusive. They still lack a comprehensive treatment, despite
the aid coming from the theory of super Harish-Chandra pairs (SHCP), which
allows to effectively %deal with a supergroup by
identify a supergroup %(categorically)
with a pair consisting of its
underlying group and its Lie superalgebra \cite{MS}, \cite{G}, \cite{CF}.

\medskip
In the ordinary setting, a reductive algebraic group $G$ over an
algebraically closed field is completely described by
its root datum \cite{Sp}. Once we have fixed a maximal torus $T \subset G$,
we can compute the characters $X(T)$ and cocharacters $X^\vee(T)$,
and then define the root datum $(X(T), \Phi(T), X^\vee(T), \Phi^\vee(T))$,
with $\Phi(T)$ denoting the root system of the semisimple group associated to $G$.
%is then naturally associated with $G$.
Vice-versa, given an axiomatic definition of
root datum, we can %show that we can
always associate to it a connected
reductive algebraic group and two groups are isomorphic if and only
if their root data are equivalent (see \cite{Sp} for more details).

\medskip
In the present paper we want to tackle two distinct and yet related
questions, that generalize to the supersetting
the bijective correspondence between
root data and connected reductive algebraic groups, briefly described above.

\medskip
  {First we want to investigate the pairs $(G,Y)$, consisting
of a connected reductive algebraic group $G$ and an odd
G-superscheme $Y$, that come from a supergroup $\bbG$  {over an algebraically closed field $\bk$ of characteristic $0$},
that is $G=\bbG_\ev$, the underlying group,
and $Y=\bbG/\bbG_\ev$, with the natural $G$ action. We also call
$(G,Y)$ the \textit{quasi-reductive pair} of $\bbG$}.

We identify a meaningful class of such pairs, which we call
\textit{basic quasi-reductive pairs},  by means of the notion of
{\it generalized root system}
as in \cite{Ser} (see also \cite{DF}, \cite{DF2} for another point of view
on root systems of Lie superalgebras).
%In fact, considering the action of a maximal torus $T \subset G$ on
%the tangent space to $Y$ at its only topological point, we require the
%weight spaces to satisfy the generalized root axioms in \cite{Ser}.
%we require the torus action, for a maximal torus $T \subset G$,
%on the tangent space to $Y$ at its only topological point, satisfies the relevant
%generalized root system definition.
  {We say that a supergroup is {\it basic quasi-reductive}, if its quasi-reductive
pair $(\bbG_\ev, \bbG/\bbG_\ev)$ is basic quasi-reductive.}
%\comment{removed ``indecomposable'', to ease the reading}
%Furthermore, we can talk about an indecomposable basic quasi-reductive pair $(G,Y)$
%(see the paragraph prior to Corollary \ref{cor: 2.9}).

\medskip
Our first main result is the following (Theorem  \ref{prop: basic existence}
and Corollary \ref{cor: 2.9}).

\begin{theorem}\label{main1}
Let $(G, Y )$ be an indecomposable  basic quasi-reductive pair.
Then, there exists an algebraic supergroup $\bbG$
with $(G,Y)$ as its quasi-reductive pair.
Furthermore, such $\bbG$ is unique,
up to isomorphisms.
\end{theorem}

Then, we turn to examine a classification question for a class of quasi-reductive supergroups {{satisfying quasi-reductive pair conditions (see Definition \ref{defn: qrp})}}. For these supergroups
we have that their root data completely characterizes them
and we can prove the following classification result (Theorems \ref{thm: 3.11},
\ref{thm: 2.13}). For the notion of isogeny class see after Prop. \ref{prop: 3.10}.

\begin{theorem}\label{main2}
\begin{enumerate}
\item Two basic quasi-reductive supergroups are isomorphic
if and only if they are associated to equivalent root data.
\item  An indecomposable  basic quasi-reductive supergroup {{belongs to one isogeny class
%is isomorphic to
  whose  representative}}  comes from  the following list:
\begin{itemize}
\item[(2.1)]  %$\text{GL}(m,n)$ with $m,n>0$;
$\text{SL}(m,n)$ with $m>n>0$; {$\widetilde{\text{SL}}(n,n)$ with   $n>2$,}
where
 $\text{SL}(m,n)$ is defined via for $A\in \salg$
 \begin{align*}
 \text{SL}(m,n)(A)=&\{g=\left( \begin{array}{cc}
B_1 & B_2\cr
B_3 & B_4
\end{array}\right)
 \in \text{GL}(m,n)(A)\mid \textsf{Ber}(g)=1 \cr
 &\;\;\text{ with } \textsf{Ber}(g):=\textsf{det}(B_1-B_2B_4^{-1}B_3)\textsf{det}
 (B_4)^{-1}\}
 \end{align*}
 and
$\widetilde{\text{SL}}(n,n)$ is defined via for $A\in \salg$
$$\widetilde{\text{SL}}(n,n)(A)=\text{SL}(n,n)(A)\slash \langle aI_{2n}\mid a\in U(A_\bz)\rangle$$
{ with $\salg$ denoting the category of commutative $\bk$-superalgebras,  and $U(A_\bz)$ denoting the subset of invertible elements in  $A_\bz$;}
\item[(2.2)] $\text{OSp}(2m+1,2n)$ with $m\geq 0, n\geq1$,
{{here and below $OSp$ indicates the corresponding connected supergroup. }}
\item[(2.3)] $\text{OSp}(2,2n)$ with $n\geq 1$;
\item[(2.4)] $\text{OSp}(2m,2n)$ with $m\geq 2$, $n\geq 1$;
\item[(2.5)] $D(2,1,\alpha)$, $\alpha$ is a nonzero scalar;
\item[(2.6)]  $F(4)$;
\item[(2.7)]  $G(3)$.
\end{itemize}
\end{enumerate}
\end{theorem}

%Finally, for a characteristic zero algebraically
%closed field, we consider the case of {\it reductive supergroups}, %$\bbG$,
%that is supergroups with trivial unipotent radical.

  {Then, in the hypothesis of a nondegenerate
bilinear form on the Lie superalgebra $\fg=\Lie(\bbG)$
of a given {quasi-reductive} supergroup $\bbG$ {along with the invertibility hypothesis  of all odd reflections}, we prove
that $\bbG$ is basic quasi-reductive of {\sl monodromy type} (Theorem \ref{thm: 4.1}).
$\bbG$ of monodromy type means %that is those groups
that the action of its maximal torus on the space $\bbG/\bbG_\ev$ has
one dimensional weight spaces.}

\begin{theorem}\label{main3}
Let $\bbG$ be a connected quasi-reductive supergroup,
$T$ a maximal torus of its underlying group $\bbG_{\mathrm{ev}}$. If %its Lie superalgebra
$\fg=\Lie(\bbG)$ admits a nondegenerate even
super-symmetric invariant bilinear form and {the invertible property of all odd reflections}, then  $\bbG$ is
a basic quasi-reductive supergroup of monodromy type.
%  {Here a
%{\it supergroup of monodromy type} means %that is those groups
%that the action of its maximal torus on the space $\bbG/\bbG_\ev$ has
%one dimensional weight spaces. }
\end{theorem}

{\bf Acknowledgements.} Both R.F. and B.S. would like to thank the
FaBiT (Farmacia e Biotecnologie) and
DIFA (Dipartimento di Fisica e Astronomia)
at the University of Bologna, for the warm hospitality
during the completion of this work.
R.F. acknowledges the funding Gnsaga-Indam, by COST Action
CaLISTA CA21109, MSCA-2022-SE CaLIGOLA 101086123, MSCA-2023-DN CaLiForNIA - 101119552,
PNRR MNESYS, PNRR SymQuSec, PNRR National Center for HPC, Big Data and
Quantum Computing, INFN Sezione Bologna. B.S. thanks An Zhang for helpful discussions, also he acknowledges support by the National Natural Science Foundation of China ({Nos. 12071136, 12271345,}), and by Science and Technology Commission of Shanghai Municipality (No. 22DZ2229014).
{ Both  R.F. and B.S. would like to thank the Referees for the
  careful reading of the manuscript.}

\section{Preliminaries}
%{Purely-odd quotient superschemes and super partners for algebraic supergroups}
Throughout,  {{$\bk$ is assumed to be an algebraically  closed field of  characteristic $0$. }}
%not equal to $2$
%{ or $3$.}
Let $\salg$ denote the category of commutative $\bk$-superalgebras.
For $A=A_\bz\oplus A_\bo\in \salg$, denote by $U(A_\bz)$ the group of invertible elements in $A_\bz$
and with $J_A$ the ideal generated by odd nilpotents.
If $X=(|X|,\cO_X)$ is an affine superscheme, we denote with $|X|$ the underlying
topological space, with $\cO_X$ its structural sheaf and
with $\cO(X)$ the superalgebra of global sections
(see \cite{CCF}, \cite{Var} for notation and preliminaries).

\subsection{Algebraic Supergroups}
%{Purely-odd quotient superschemes and the global splitting for algebraic supergroups}
Let $\bbG$ be an algebraic supergroup over $\bk$ with Lie superalgebra $\ggg=\ggg_\bz\oplus\ggg_\bo$, and purely-even subgroup $G:=\bbG_\ev$, where $\bbG_\ev$ is the ordinary algebraic group underlying $\bbG$, i.e. $\cO(\bbG_\ev)=\cO(\bbG)/J_{\cO(\bbG)}$.
Then $\ggg_\bz=\Lie(G)$.

Note that $G$ is a closed supersubgroup scheme of $\bbG$. We have the
quotient superscheme $$Y:=\bbG\slash G$$
(see \cite{MZ}, and \cite{MT}, \cite{FKT}, \cite{Br}, and refs. therein),
which is a purely-odd superscheme %\blue{with  $G$-action},
i.e. $|Y|=\{y\}$ is just a point,   {carrying a natural $G$-action}.
We have the following result.

\begin{theorem}\label{thm: Gavarini} (\cite[Theorem 3.2.8, Corollary 3.2.9, Proposition 4.2.13, Proposition 4.3.5]{G})
\begin{itemize}
\item[(1)] There is an isomorphism of $\bk$-superschemes: $$\bbG\cong G\times Y.$$

\item[(2)] As a superscheme, $Y$ is isomorphic to $\bba_\bk^{(0|N)}$ for $N=\dim \ggg_{\bo}$.
   More precisely,
    %\blue{let $\bbG\subset
%    $Y$ can be described when $\bbG$ is realized as a closed subgroup of the general linear supergroup
%$\bbG$ can be regarded a closed subgroup of $\bbGL(V)$  for  $V=\textsf{ind}^\ggg_{\ggg_\bz}\textbf{1}\cong\bigwedge^\bullet(\ggg_\bo)$
 %   where $\textbf{1}$  is a one-dimensional trivial $\ggg_\bz$-module. Then,
 for any $A\in \salg$
    $$Y(A)={\overset\rightharpoonup{\prod}}_{i\in I}(1+A_\bo\Theta_i)=\{{\overset\rightharpoonup{\prod}}_{i\in I}(1+a_i\Theta_i)\mid a_i\in A_\bo, \forall i\in I\}  $$
 where $\Theta_i$ ($i\in I$) is a given basis of $\ggg_\bo$ with the ordered index set $I$, and ${\overset\rightharpoonup{\prod}}_{i\in I}  $
denotes an ordered product.

    \item[(3)] The global sections $\co_Y = \bigwedge^{\bullet}\ggg_\bo$.
    \end{itemize}
\end{theorem}
From the above result, it follows that  $Y$ is actually a super subscheme of $\bbG$.

\begin{remark}
According to the super Harish-Chandra pair theory, there is a category equivalence
between the categories of algebraic supergroups and of
super Harish-Chandra pairs
(\cite[\S7.4]{FKT}, \cite[\S4.3]{G}, \cite{M}, \cite{CF} and \cite[Theorem 4.22]{MS}).
%, {\sl{etc}}.).
\end{remark}

Let $\bbG$ be an algebraic supergroup scheme over $\bk$, $Z$ be a superscheme over $\bk$.
%Let $\salg$ denote the category of commutative $\bk$-superalgebras.

\begin{defn}\label{defn: action}
We say that $\bbG$ acts on $Z$ if we have a morphism
$\tau: \bbG \times Z \rightarrow Z$ corresponding to the functorial family of morphisms
$$\tau_A: \bbG(A)\times Z(A)\rightarrow Z(A), (g,z)\mapsto g.z \;\;\forall g\in \bbG(A), z\in Z(A) \text{ for  any } A\in \salg$$
satisfying the following properties:
\begin{itemize}
\item[(1)] $e.z=z$ for any $z\in Z(A)$ where $e$ is the identity of $\bbG(A)$.
\item[(2)] $(g_1g_2).z=g_1.(g_2.z)$ for any $g_1,g_2\in \bbG(A)$ and $z\in Z(A)$.
\end{itemize}
\end{defn}

\section{Basic quasi-reductive root data for  basic quasi-reductive supergroups}
In the following we always suppose $G$ is a
connected reductive algebraic group over $\bk$.

\subsection{Notation.}

\label{sec: basic notation}

Let $G$ be a connected reductive group over $\bk$ with $\ggg=\Lie(G)$,
let $T$ be a fixed maximal torus of $G$, with related root system
$\Phi:=\Phi(T)$, character group $X(T)$ and cocharacter group $X^\vee(T)$. Recall that there is a duality pairing
$$X(T)\times X^\vee(T)\rightarrow \bbz,\;(\lambda,\chi)\mapsto \langle\lambda,\chi\rangle$$
such that $\lambda(\chi(t))=t^{\langle\lambda,\chi\rangle}$ for any $\lambda\in X(T)$, $\chi\in X^\vee(T)$ and $t\in \bk^\times:=\bk\backslash\{0\}$. This pairing gives rise to a bijection of $X(T)$ to $X^\vee(T)$, $\lambda\mapsto \lambda^\vee$.

Let $Y$ be a purely-odd superscheme of super-dimension $(0|N)$,
$|Y|=\{y\}$. Suppose $Y$ is endowed with $G$-action such that
$\epsilon:\mathcal{O}(Y) \longrightarrow k$,
 $\epsilon:=\ev_y$, evaluation at $y$
is a $G$-equivariant map. {Its cotangent space is defined as (\cite{CCF}, Ch. 5):
\begin{align}\label{eq: cotangent}
  %T^*_y(Y)
  D(Y):={\mathrm{ker}}(\epsilon)\slash ({\mathrm{ker}}(\epsilon))^2
\end{align}
It is a  finite dimensional
superspace},  which will be simply denoted by  $D$
{ whenever there is no danger of confusion.}
%or also by $D(Y)$ to mark its dependence on $Y$}.
Then $D$  becomes a rational $G$-module.
The following lemma is clear.

\begin{lemma}\label{lem: odd weights} As a $G$-module,
$D$ can be decomposed into a direct sum of  weight spaces
\beq\label{D-dec}
D=\sum_{\gamma\in \Gamma(D)\subset X(T)}D_\gamma
\eeq
{{where
$\Gamma(D)=\{\gamma\in X(T)\mid  D_\gamma\ne0\}$.
}}
 \end{lemma}

%\green{The weights from the above $\Gamma:=\Gamma(D)$ will be  called
%\textit{purely-structural}  weights.}
%\commento{R. I would remove, we never use this terminology later}

%\blue{
\begin{definition}\label{quasi-red-def}
We say that $\bbG$ is a \textit{quasi-reductive supergroup}
if its underlying ordinary group $\bbG_\ev$ is reductive.
\end{definition}

If $\bbG$ is a quasi-reductive supergroup, with $G=\bbG_\ev$,
%underlying ordinary group and $Y=\bbG/\bbG_\ev$. T
then by Theorem
\ref{thm: Gavarini}, we have that $Y=\bbG/\bbG_\ev$
is a purely odd superscheme,
carrying a natural $G$-action.

Based on our previous treatment,
we give the following definition, that we shall use in the sequel.

\begin{definition}\label{root-sys-def}
Let $\bbG$ be a quasi-reductive supergroup and
$G=\bbG_\ev$, $Y=\bbG/\bbG_\ev$. Let $T$ be a maximal torus of $G$.
We define, with an abuse of notation, the \textit{root system} of $\bbG$ as
$\mathfrak{R}=\Phi(T) \cup \Gamma(D)$, where $\Phi(T)$ is the root system
of $G$ and $\Gamma(D)$ is defined in (\ref{D-dec}).

We shall write $\Phi$ and $\Gamma$ for short, whenever there is
no danger of confusion.
\end{definition}

We will see  that $\mathfrak{R}$ shares some properties of
root systems of complex semi-simple Lie algebras;
we are going to describe such properties more in details in the next sections,
adding some extra hypotheses.

\subsection{Basic quasi-reductive pairs}

We give a definition, which is the key for our subsequent treatment.

\begin{defn}\label{defn: qrp} Let $G$ be a connected reductive algebraic group,
and $Y$ {{an affine}} $G$-superscheme of super-dimension $(0|N)$.
Call the pair $(G,Y)$ a \textit{basic quasi-reductive pair}
if it satisfies the following axioms:
\begin{itemize}
\item[BQR(1)] $\Gamma=-\Gamma\subset X(T)\backslash \{0\}$,
{where {{$\Gamma=\Gamma(D(Y))$}} is the set of weights for the $G$-module $D(Y)$, see (\ref{eq: cotangent}).}
%$T^*_y(Y)$ for $y\in |Y|$ in the same sense of  (\ref{eq: cotangent}).}
\item[BQR(2)] Set $\frakr:=\Phi\cup\Gamma$, then %The finite set
$\frakr$ spans $X(T)$.
\item[BQR(3)] Denote by $V$ and $V_0$  the $\bbr$-space spanned by $\frakr$ and $\Phi$, respectively. There is a nondegenerate symmetric  bilinear form $(\cdot,\cdot)$ on $V$ such that
 $(\cdot,\cdot)|_{V_0}$ is compatible with the duality  pairing $\langle\cdot,\cdot\rangle$: $X(T)\times X^\vee(T)\rightarrow \bbz$, which means that $\langle\alpha,\beta^\vee\rangle=2(\beta,\beta)^{-1}(\alpha,\beta)$ for any $\alpha,\beta\in \Phi$.
\item[BQR(4)]
%(Serganova's condition)
For $\gamma\in \Gamma$,
\begin{itemize}
\item[(S1)]  if $(\gamma,\gamma)=0$,  there is an invertible mapping $$\sfr_\gamma: \frakr\rightarrow \frakr$$
 such that for $\nu\in \frakr$
 \begin{align*}
 \sfr_\gamma(\nu)=\begin{cases} \nu\pm \gamma &\text{ if } (\gamma,\nu)\ne 0;\cr
   \nu &\text{ if } (\gamma,\nu)= 0.
  \end{cases}
\end{align*}
\item[(S2)] If $(\gamma,\gamma)\ne0$, then $\langle \nu,\gamma^\vee\rangle= 2(\gamma,\gamma)^{-1}(\nu,\gamma)$ for $\nu\in\frakr$, and
$\nu-\langle \nu,\gamma^\vee\rangle \gamma\in \frakr$.
%; and $\nu^\vee-\langle \gamma,\nu^\vee\rangle \gamma^\vee\in %\frakr^\vee$
\end{itemize}
\end{itemize}
\end{defn}
Compare (BQR4) with Def. 1.1 by Serganova in \cite{Ser1}.

%\subsubsection{Examples}

\begin{observation}  Let $\bbG$ be a
basic classical supergroup, that is  %with $G:=\bbG_\ev$,
$\ggg:=\Lie(\bbG)$ is a basic classical superalgebra { which is not  $A(1|1)$}
(see \cite{Kac}, \cite{FG}, and refs. therein).
Keeping in mind Theorem \ref{thm: Gavarini}(3) along with the properties of root systems
of basic classical Lie superalgebras (see \cite{Kac}, \cite{CW}, \cite{M} and refs. therein),
we see that $(G, Y:=\bbG\slash G)$ becomes a basic quasi-reductive pair,
{with $G:=\bbG_\ev$.}
{Moreover, in the light of Thm. 7.2 in \cite{Ser1} and the
classification in \cite{FG} we will see that non isomorphic basic classical supergroups will
give different, i.e. non isomorphic, quasi-reductive pairs.
We are going to see more precisely this statement and its converse in Thm \ref{prop: basic existence}.}%thm: 3.11}.}
\end{observation}

\subsection{Basic quasi-reductive root data}\label{sec: basic qrrd}

We want to associate to a basic quasi-reductive pair a {\sl root data}, in analogy to
theory of reductive groups (see \cite{Sp}, \cite[Ch. II.1]{Jan}).

Suppose  $(X,R, X^\vee,  R^\vee)$ is the root datum for a
connected reductive algebraic group as in \cite[\S7.4]{Sp} (or see \cite[Appendix C]{Mi}).

\begin{definition}
Let $\Gamma\subset X$ be a finite subset without intersection with $R$.
Let $\frakr=R\cup \Gamma$, and $\frakr^\vee=R^\vee\cup \Gamma^\vee\subset X^\vee$
where $R^\vee$ and $\Gamma^\vee$ are the duals of $R$ and $\Gamma$ in the pairing
of $X$ and $X^\vee$.
If $\frakr$ satisfies the axioms of BQR(1)-(4)
in Definition \ref{defn: qrp} with $R$ taking the place of  $\Phi$
and $V=\bbr\otimes_\bbz \bbz\frakr$,
then we call $(X,\frakr,X^\vee,\frakr^\vee)$  a
\textit{basic quasi-reductive datum}.
  {For simplicity we also  call $\frakr$ a \textit{basic quasi-reductive root system}.}
%  Correspondingly, $V=\bbr\otimes_\bbz \bbz\frakr$ where $\bbz\frakr$ stands for a $\bbz$-free module spanned by $\frakr$.
\end{definition}

Recall Serganova introduced generalized root systems (\GRS\footnote{
{ %\begin{remark}
We prefer the acronym \GRS i.e. Serganova Generalized Root Systems, to distinguish
them from the generalized root systems GRS as in \cite{DF2}, which are more general.
%\end{remark}
}} for short), see \cite[\S1]{Ser1}.
A \GRS is called \textit{irreducible} if it can not be decomposed into a sum of two \GRS (see \cite[Definition 1.13]{Ser1}).
She classified there all irreducible SGRS, and  proved that  the irreducible \GRS are exactly the root systems of basic classical  Lie superalgebras.

\begin{lemma}\label{lem: 3.4} Retain the above notations. The following the statements hold.
\begin{itemize}
\item[(1)] The root system  $\frakr\subset V\backslash\{0\}$ forms a generalized root system in the sense of Serganova \cite{Ser1}.

\item[(2)] Furthermore, $\frakr$ can be decomposed into a direct sum of  irreducible generalized root systems in  Serganova's sense \cite{Ser1}, i.e. a union of irreducible \GRS without intersection for any two different irreducible components.
\end{itemize}
\end{lemma}

\begin{proof} (1) Under Serganova's condition BQR(4) in Definition \ref{defn: qrp},  it can be checked by using the property of reductive algebraic groups.

(2) Note that \cite[Lemma 1.14]{Ser1} is suitable to ordinary cases, independent of the characteristic of the base field. This part follows from (1)  and \cite[Lemma 1.14]{Ser1}.
\end{proof}

\begin{observation} Suppose $(G,Y)$ is a basic quasi-reductive pair of the basic quasi-reductive  supergroup $\bbG$. Maintain notations as in \S\ref{sec: basic notation}. Then
$$\Psi(G,Y,T):=(X(T), \frakr(T,D), X^\vee(T), \frakr(D,T)^\vee)$$ becomes
a basic quasi-reductive root datum, where $\frakr(D,T)=\Phi(G,T)\cup \Gamma(D)$.
\end{observation}

  {%By Lemma \ref{lem: 3.4},
We say that a basic quasi-reductive datum
$(X,\frakr,X^\vee,\frakr^\vee)$ is \textit{irreducible}
if $\frakr$ is irreducible, in the sense of \cite{Ser1}.
We say that a basic quasi-reductive pair $(G,Y)$ is \textit{irreducible}
if its root datum $\frakr$ is irreducible.

\medskip
These notions are well defined
because of Lemma \ref{lem: 3.4}.
}
%\subsection{Equivalence}

\begin{defn}\label{def:eq}
We say that two basic quasi-reductive data
$(X_i,\frakr_i, X_i^\vee,\frakr_i^\vee)$, $i=1,2$,  are \textit{equivalent} if there is an isomorphism
{ $\phi:X_1 \lra X_2$}
mapping $R_1$ onto $R_2$, $\Gamma_1$ onto $\Gamma_2$ and such that the dual morphisms
map $R_2^\vee$ onto $R_1^\vee$, and $\Gamma_2^\vee$ onto $\Gamma_1^\vee$.
{{For two basic quasi-reductive data, we further}} have the commutative diagram:
\beq
\xymatrix{
    \frakr_1 \ar[r]^\phi \ar[d]_{\sfr_{\gamma_1}} & \frakr_2 \ar[d]^{\sfr_{\gamma_2}} \\
    \frakr_1 \ar[r]_{{\phi}}       & \frakr_2 }
%\sfr_\gamma
\eeq
for ${\gamma_i} \in \frakr_i$.
\end{defn}

% We say $R$ is a direct sum of $R'$ and $R''$ for  finite set $R'$ and %$R''$ in $X(T)$ if $V=V'\oplus V''$

\subsection{Basic quasi-reductive supergroups}

We show how basic reductive root data is encoding all the information
regarding a supergroup
belonging to a class we call {\sl basic quasi-reductive supergroups}.

%\subsection{Basic quasi-reductive supergroups}
\begin{defn}\label{defn: 3.3}  Let $\bbG$ be an algebraic supergroup over $\bk$.
We call $\bbG$ a \textit{basic quasi-reductive supergroup}
if  $\bbG_\ev$ is a connected algebraic group, and  $(\bbG_\ev,Y)$
is a basic quasi-reductive pair for $Y=\bbG\slash \bbG_\ev$.
\end{defn}

%From now on till the end of this section, we assume that $\text{ch}(\bk)=0$.

\begin{theorem}\label{prop: basic existence}
    {(Existence and Uniqueness Theorem).
Let $(G,Y)$ be an irreducible basic quasi-reductive pair.
Then, there is an algebraic supergroup  $\bbG$ with $(G,Y)$ as its quasi-reductive pair.
%$\bbG_{\mathrm{ev}}=G$ such that $Y=\bbG/\bbG_{\mathrm{ev}}$.
%is a supergroup partner of $G$.
Furthermore, such  $\bbG$ is unique, up to isomorphisms.}
\end{theorem}

\begin{proof} According to Lemma \ref{lem: 3.4}, $\frakr$ is a generalized root system in the sense of Serganova \cite{Ser1}. Hence, by applying %Serganova's result
\cite[Theorem 7.2]{Ser1}
along with the same arguments as in %Fioresi-Gavarini's paper
\cite[Chapters 5 and 6]{FG}, using the Chevalley basis and the integral
Kostant form the existence is proved over an arbitrary field.
As to the uniqueness, it follows
from the theory of Super Harish-Chandra Pairs (SHCP). %Lemma \ref{lem: gen iso}.
\end{proof}

Notice that $\bbG$ in Theorem \ref{prop: basic existence} is by definition basic quasi-reductive.

\medskip

%\begin{defn}
%Such an algebraic supergroup $\bbG$ which is showed to exist in  Theorem \ref{prop: basic existence} is  called a basic quasi-reductive algebraic supergroup.
%\end{defn}
%\subsection{}\label{sec: 2.6}

We say that $(G,Y)$ is \textit{decomposable} if there are  two  connected reductive closed subgroups
$G'$ and $G''$ of $G$ with $G= G'\times G''$, and two {{ $G$-subsuperschemes
$Y'$ and $Y''$ of $Y$ with $Y\cong Y'\times Y''$ as $G$-superschemes}} such that $(G', Y')$ and $(G'', Y'')$
become basic quasi-reductive pairs. We call $(G,Y)$ \textit{indecomposable} if otherwise.
%By Theorem \ref{prop: basic existence}, since a basic quasi-reductive pair is associated to a unique
%basic quasi-reductive supergroup, we say that $
We say that a basic quasi-reductive supergroup $\bbG$ is  \textit{indecomposable} if
$(\bbG_\ev, \bbG/\bbG_\ev)$ is indecomposable.

\begin{corollary}\label{cor: 2.9}
Let $(G,Y)$ be a basic quasi-reductive pair.
Suppose  $(G,Y)$ is indecomposable, then $\frakr$ is irreducible.
\end{corollary}

\begin{proof} We only need to show that if $\frakr$ is decomposed into a direct sum of two irreducible \GRS $\frakr_1$ and $\frakr_2$, then $(G,Y)$ admits the corresponding decomposition.
%By the Existence and Uniqueness Theorem \ref{prop: basic existence},
By applying Serganova's result \cite[Theorem 7.2]{Ser1} along with
%Fioresi-Gavarini's paper
\cite[Chapters 5 and 6]{FG},
there are supergroups $\bbG_1$ and $\bbG_2$ such that $G_i=(\bbG_i)_\ev$ are connected reductive algebraic groups with maximal torus $T_i$, $i=1,2$, and the root systems of $\bbG_i$ associated with $T_i$  are  $\frakr_1$, and $\frakr_2$, respectively.  The root system of $\frakr$ decomposing into a direct sum of $\frakr_1$ and $\frakr_2$ implies that $\tilde T:=T_1\times T_2$  is a maximal torus of $\tilde G:=G_1\times G_2$.  By a direct routine computation on affine coordinate superalgebras it is readily seen that  there is an isomorphism of affine algebraic groups
\begin{align}\label{eq: even iso}
\tilde G\cong \tilde\bbG_\ev
\end{align}
where $\tilde \bbG:=\bbG_1\times \bbG_2$,  and  $X(T_i)$ is spanned by $\frakr_i$, $i=1,2$. Furthermore, there is an isomorphism of superschemes
\begin{align}\label{eq: Y decomp}
\tilde Y:=\tilde \bbG\slash \tilde G\cong \bbG_1\slash G_1\times \bbG_2\slash G_2=:Y_1\times Y_2
\end{align}
which follows from  Theorem \ref{thm: Gavarini}(2).
Then according to the structure theorem of connected reductive algebraic groups (see \cite[Theorem 7.3.1 and Corollary 7.6.4]{Sp}), it is readily concluded that $G$ is isomorphic to $(G_1\times C_1)\times (G_2\times C_2$), where $C_i$ are some central tori of $G$ respectively. Combining with (\ref{eq: Y decomp}), we obtain a contradiction with the indecomposable assumption of $(G,Y)$,
  {taking into account Thm \ref{prop: basic existence}}.
  The proof is completed.
%(2) For the sufficient part,  we only show that if $(G,Y)$ is  decomposable, %then $\frakr$ is not irreducible. This follows from Lemma \ref{lem: 3.4}.
\end{proof}

  {
On the other hand, we defined the notion of
irreducible quasi-reductive root datum $$\Psi:=(X,\frakr,X^\vee,\frakr^\vee),$$
(see Lemma \ref{lem: 3.4}(1)), by asking
$\frakr$ to be an irreducible SGRS.
By Serganova's result in \cite{Ser1} and the same arguments as in
\cite[Chapters 5 and 6]{FG},
there is  an algebraic supergroup $\bbG$ such that $\bbG_\ev=G$,
and  $(\bbG_\ev, \bbG\slash\bbG_\ev)$ is a quasi-reductive pair
associated to which $\Psi$ is the root datum. Hence,
by Theorem \ref{prop: basic existence} we have the following.}

\begin{prop} \label{prop: 3.10}
For any irreducible quasi-reductive root datum
$\Psi=(X,\frakr,X^\vee,\frakr^\vee)$,
there is a unique (up to isomorphisms) algebraic supergroup $\bbG$
whose quasi-reductive root datum is $\Psi$.
\end{prop}

  {
There are possibly different irreducible quasi-reductive root data
$\Psi_1$ and $\Psi_2$ associated with the same root system $\frakr$.
So there are possibly different (up to isomorphisms) indecomposable
algebraic supergroup $\bbG_1$ and $\bbG_2$ associated with the same root
system $\frakr$.  All  such $\bbG_i$ associated with $\frakr$ form a
class of indecomposable  basic quasi-reductive algebraic supergroups,
called {\textit{ the isogeny class of $\frakr$}}.
}

%\subsection{The isomorphism theorem}
%\subsection{Equivalence}
%\begin{defn} We say that two basic quasi-reductive data  $(X_i,\frakr_i, X_i^\vee,\frakr_i^\vee)$  are equivalent if there is an isomorphism $X_1$ onto $X_2$ mapping $R_1$ onto $R_2$, $\Gamma_1$ onto $\Gamma_2$ such that the dual maps $R_2^\vee$ onto $R_1^\vee$, and maps $\Gamma_2^\vee$ onto $\Gamma_1^\vee$.
%\end{defn}

In Def. \ref{def:eq} we established the notion of equivalence of
two basic quasi-reductive data  $(X_i,\frakr_i, X_i^\vee,\frakr_i^\vee)$. Next theorem states
that equivalence of root data correspond to isomorphism of supergroups, in the appropriate hypotheses.

\begin{theorem}\label{thm: 3.11}
Let $\bbG_i$ ($i=1,2$) be two basic quasi-reductive supergroups with purely-even groups $G_i$, $i=1,2$ respectively. The basic quasi-reductive data $\Psi(G_i,Y_i,T_i)$ are equivalent,
if and only if $\bbG_i$, $i=1,2$ are isomorphic as algebraic supergroups.
\end{theorem}

\begin{proof} By the decomposable property of generalized root systems (see \cite[Lemma 1.14]{Ser1}), it can be reduced to the case of indecomposable basic quasi-reductive supergroups. Then by Proposition \ref{prop: 3.10}, we can accomplish the remaining arguments for the proof of the  theorem.
The other direction is clear.
\end{proof}

\subsection{Classification of indecomposable basic quasi-reductive supergroups}

%\vskip5pt
%\subsubsection{}
Thanks to Serganova's classification on \GRS in \cite{Ser1},
we have the following  classification theorem, by summarising the above results.

\begin{theorem}\label{thm: 2.13}
Let $\bbG$ be an indecomposable  basic quasi-reductive supergroups over $\bk$ {{which is assumed to be an algebraically closed field of characteristic $0$.}}
Then $\bbG$ is isomorphic to a supergroup
whose isogeny class has a representative in the following list:
\begin{itemize}
\item[(2.1)]  $\text{SL}(m,n)${{ with $m> n>0$,}}  %$\widetilde{\text{SL}(n,n)}$ with   $n\ne2$,
%$\text{GL}(m,n)$; $\text{SL}(m,n)$ with $m,n>0$; %$\widetilde{\text{SL}(n,n)}$ with   $n>0$,
where
 $\text{SL}(m,n)$ is defined via for $A\in \salg$
 \begin{align*}
 \text{SL}(m,n)(A)=&\{g=\left( \begin{array}{cc}
B_1 & B_2\cr
B_3 & B_4
\end{array}\right)
 \in \text{GL}(m,n)(A)\mid \textsf{Ber}(g)=1 \cr
 &\;\;\text{ with } \textsf{Ber}(g):=\textsf{det}(B_1-B_2B_4^{-1}B_3)\textsf{det}
 (B_4)^{-1}\}
 \end{align*}
 and
$\widetilde{\text{SL}}(n,n)$ {{with $n>2$}}, which is defined via for $A\in \salg$
$$\widetilde{\text{SL}}(n,n)(A)=\text{SL}(n,n)(A)\slash \langle aI_{2n}\mid a\in U(A_\bz)\rangle.$$
{(Recall that $\salg$ denotes the category of commutative $\bk$-superalgebras,  and $U(A_\bz)$ denotes the subset of invertible elements in  $A_\bz$.)}
\item[(2.2)] $\text{OSp}(2m+1,2n)$ with $m\geq 0, n\geq1$,
{{here and below $OSp$ indicates the corresponding connected supergroup. }}
\item[(2.2)] $\text{OSp}(2m+1,2n)$ with $m\geq 0, n\geq1$
\item[(2.3)] $\text{OSp}(2,2n)$ with $n\geq 1$
\item[(2.4)] $\text{OSp}(2m,2n)$ with $m\geq 2$, $n\geq 1$
\item[(2.5)] $D(2,1,\alpha)$, $\alpha$ is a nonzero scalar.
\item[(2.6)]  $F(4)$
\item[(2.7)]  $G(3)$
\end{itemize}
\end{theorem}

\section{Basic quasi-reductive supergroup of monodromy type}

%In this section, the characteristic of $\bk$ is assumed to be zero.
  {We consider the class of %$\bbz$-splitting
quasi-reductive algebraic supergroups with
a nondegenerate even super-symmetric
invariant bilinear form \cite{Kac} on their Lie superalgebras.}

\subsection{Quasi-reductive supergroups with nondegenerate bilinear forms}

%In the main result of this section, Theorem \ref{thm: 4.1},
%we characterize connected
%quasi-reductive supergroups (see Def. \ref{quasi-red-def})
%with a nondegenerate super-symmetric invariant form.

%  {as basic quasi-reductive supergroup of monodromy type.

Recall that for  a quasi-reductive supergroup $\bbG$, we can naturally
associate to it the root system $\mathfrak{R}=\Phi(T) \cup \Gamma(D)$,
as in Def. \ref{root-sys-def}, where $\Phi(T)$ is the root system
of the purely-even subgroup group $G$ ($T$ maximal torus of $G$)
and $\Gamma(D)$ is defined in
(\ref{D-dec}). We also write $\Phi(G,T)$, or simply $\Phi$, for $\Phi(T)$.
With an abuse of terminology,
we call the elements in $\Phi(T)$ \textit{even roots} and those in $\Gamma(D)$
\textit{odd roots}   {and we call $\Phi(T) \cup \Gamma(D)$ the
\textit{root system} of $\bbG$.}

\iffalse
 \blue{
\begin{remark}We remind the reader that a connected quasi-reductive supergroup does not contain any component one isomorphic to the supergroup of $A(1|1)$.
 \end{remark}
}
\fi

%The proof of Theorem \ref{thm: 4.1} will be given in the following  subsections.
%We give a few preliminaries.
 % in the following subsections.
%We first establish some notation.

%\medskip
 Let $\ggg$ be the Lie superalgebra of $\bbG$.
 Consequently, $\ggg=\ggg_\bz\oplus \ggg_\bo$ with $\ggg_\bz\cong \Lie(G)$. By assumption,  $G$ is a connected reductive algebraic group, and $\ggg_\bz$ a reductive Lie algebra. Furthermore, suppose $\ttt=\Lie(T)$.
Simply write $\Phi=\Phi(G,T)$ and $\Gamma=\Gamma(D)$ as before. Set $\scrw$ to be the abstract Weyl group of $G$. Note that associated with $T$, $G$ has the root datum $(X(T), \Phi, X^\vee(T), \Phi^\vee)$ where $X^\vee(T)$ denotes the co-character group of $T$ (see \cite[\S7.3.4]{Sp}), and $\Phi^\vee:=\Phi^\vee(G,T)$ denotes the coroots of $G$ associated with $T$ (see \cite[\S7.1.8]{Sp}). By the classical theory \cite{Sp},
there is a good pairing $\langle \cdot,\cdot \rangle
$ between $X(T)$ and $X^\vee(T)$ satisfying the axioms of the root data as in \S\ref{sec: basic notation}. In particular, the pairing gives rise to the bijection $\lambda\mapsto \lambda^\vee$ from $X(T)$ onto $X^\vee(T)$, sending the root system of even roots
$\Phi$ onto $\Phi^\vee$.

%\subsection{}
By assumption,  we can write
\begin{align}\label{eq: multiplicity ques}
\ggg=\ttt\oplus \sum_{\alpha\in \Phi}\ggg_\alpha\oplus \sum_{\gamma\in \Gamma}\ggg_\gamma.
\end{align}
In particular $\ggg_\bo=\sum_{\gamma\in \Gamma}\ggg_\gamma$. {From now on, we assume
\begin{align}\label{eq: nondeg assump}
\qquad \ggg \text{  admits a nondegenerate even
super-symmetric } \ggg\text{-invariant bilinear form}.
\end{align}
}

%By the theory of super Harish-Chandra pairs, $\bbG$ can be realised as a %closed subsupergroup of $\mathbb{GL}(\textsf{V})$ where  %$\textsf{V}=\textsf{ind}^\ggg_{\ggg_\bz}\textbf{1}
%\cong\bigwedge^\bullet(\ggg_\bo)$ where $\textbf{1}$  is a %one-dimensional trivial $\ggg_\bz$-module (see Theorem \ref{thm: %Gavarini}, \cite[Proposition 4.2.13]{G}).

By the assumption of nondegenerate even super-symmetric bilinear form
$(\cdot,\cdot)$ on $\ggg$,
we  have the following fact, which is standard, but we include for completeness.
  {Let $(\cdot,\cdot)_{\ttt}$ denote the restriction of $(\cdot,\cdot)$ to $\ttt$.

\begin{lemma}
Let the notation be as above. Then under the assumption (\ref{eq: nondeg assump}),    $(\cdot,\cdot)_{\ttt}$ is nondegenerate.
%The restriction of the bilinear form to $\ttt$ must be nondegenerate, which will be denoted by $(\cdot,\cdot)_{\ttt}$.
\end{lemma}
}

    \begin{proof} From the even property of the nondegenerate bilinear form $(\cdot,\cdot)$, it follows that the restriction of this form to $\ggg_\bz$ is nondegenerate. Now $(H, \ggg_\bz)\ne 0$ for any nonzero $H\in\ttt$. So the claim can be deduced from the following
\begin{align}\label{eq: t nondeg}
(H, \ggg_\alpha)=0,\;\; \forall \alpha\in\Phi.
\end{align}
{{Note that there is  $H_\alpha\in \ttt$ such that $\alpha(H_\alpha)\ne 0$.  If $(H, X_\alpha)\ne 0$ for some $X_\alpha\in\ggg_\alpha$,  then we have
$(H, [H_\alpha, X_\alpha])\ne 0$.}}  On the other side,
$$(H, [H_\alpha, X_\alpha])=([H,H_\alpha], X_\alpha)=0$$
which is a contradiction. Hence (\ref{eq: t nondeg}) is true, which implies our claim.
\end{proof}

%\medskip\hrule\medskip
%\blue{RITA: this section needs to be fixed or removed.
%If we remove it, then we need to keep Lemma 4.7 somewhere ($R=-R$).}
%\medskip\hrule\medskip

\subsection{The root system of a quasi-reductive supergroup with nondegenerate form.}
We start with an auxiliary lemma establishing that $\frakr=-\frakr$.
%\subsubsection{Lemma for Theorem \ref{thm: 4.1}}

\begin{lemma}\label{lem: R equal -R}
Let the notation and assumption  be as above.
The root system $\frakr$ coincides with $-\frakr$.
% All root spaces in (\ref{eq: multiplicity ques}) are one-dimensional.
\end{lemma}
\begin{proof} %Note that $G=\bbG_\ev$ is a connected reductive algebraic group with root system $\Phi$. Hence all root spaces $\ggg_\alpha$ for $\alpha\in \Phi$ are one-dimensional. So it suffices to prove the statement for $\ggg_\gamma$ with $\gamma\in \Gamma$.  %By Lemma \ref{lem: %4.3}
%%%%We consider $\bbG_\gamma$ which is a supsupergroup generated by %$U_{\pm\gamma}$ and the one-dimensional torus %$T_\gamma=\gamma^\vee(\bk^\times)$. Suppose $\bbB_\gamma$ stands for the %Borel subgroup of $\bbG$ containing $T_\gamma$. Note that $\gamma^\vee: %\textsf{G}_{\text{m}}\rightarrow T_\gamma$ is an isomorphism, where %$\textsf{G}_{\text{m}}$ is  one-dimensional multiplication group.
%We proceed with arguments by steps.

%(1) For the first part of the Lemma, note that
Since $G$ is a connected reductive algebraic group,
we have $\Phi=-\Phi$. We only need to show that $\Gamma=-\Gamma$.
Keep in mind the assumption that there is a nondegenerate even super-symmetric
$\ggg$-invariant bilinear form $(\cdot,\cdot)$ on $\ggg$.
Note that $\Gamma\subset X(T)\backslash\{0\}$. For any $\gamma\in \Gamma$,
there exists $S\in \ttt$ such $\gamma(S)\ne 0$, we fix one and denote it by $S_\gamma$.
 %Furthermore, by the property of reductive Lie algebras we might as well %furthermore assume
%\begin{align}\label{eq: H gamma}
%\gamma(S)=(S_\gamma,
%S)\;\; \forall S\in \ttt.
%\end{align}
So for any nonzero $X_\gamma\in \ggg_\gamma$, we have $([S_\gamma,X_{\gamma}],\ggg_\bo)\ne 0$. Hence
\begin{align}\label{eq: gamma and minus-1}
(S_\gamma, [X_\gamma, \ggg_\bo])\ne 0.
\end{align}
Now $[X_\gamma, \sum_{\nu\in \Gamma\backslash\{-\gamma\}}\ggg_\nu]\in \ggg_\bz$.
Hence we have
%According to the property of reductive Lie algebras,
\begin{align}\label{eq: gamma and minus-2}
(S_\gamma, [X_\gamma, \sum_{\nu\in \Gamma\backslash\{-\gamma\}}\ggg_\nu])=0
\end{align}
because  $(\ttt,\ggg_\alpha)=0$ for any $\alpha\in\Phi$
(the arguments for this are the same as in  (\ref{eq: t nondeg})).
Combining (\ref{eq: gamma and minus-1}) with
(\ref{eq: gamma and minus-2}), we have
$\ggg_{-\gamma}\ne 0$.  So $-\gamma\in \Gamma$.
We are done.
\end{proof}

We now make a key observation, that we shall use in the sequel.

%\subsubsection{}
\begin{observation}
Notice that, by the above arguments, if $X_\gamma \in \ggg_\gamma
\subset \ggg_\bo$, there is an $X_{-\gamma}$ such that
$(X_\gamma,X_{-\gamma})=1$ and $[X_\gamma,X_{-\gamma}]\in \ttt$.
On the other hand, by the non-degeneracy of $(\cdot,\cdot)_{\ttt}$ there is a unique
$H_\gamma\in \ttt$ such that $\gamma(H)=(H_\gamma, H)$
for any $H\in \ttt$. Then
$$
([X_\gamma,X_{-\gamma}], H)=(X_\gamma, [X_{-\gamma},H])=\gamma(H)(X_\gamma,X_{-\gamma}).
$$
%By (\ref{eq: H gamma}) and
By the non-degeneracy of the restriction $(\cdot,\cdot)_{\ttt}$  again, we have
\begin{align}\label{eq: gamma Lie bracket}
[X_\gamma,X_{-\gamma}]=(X_\gamma,X_{-\gamma})H_\gamma=H_\gamma.
\end{align}

%\subsubsection{}
We further assert that $[\ggg_\gamma,\ggg_{-\gamma}]$ is spanned by $H_\gamma$ for any $\gamma\in\Gamma$, as precisely described in (\ref{eq: H theta}).
This can be verified by the same arguments as in semi-simple Lie theory (for example \cite[\S8.3]{Hum1}).
%\commento{R. I would not refer to a later eq.}
%. For the author's convenience, we make it into an account. For any %$X\in \ggg_\gamma$, $Y\in\ggg_{-\gamma}$ and $H\in\ttt$, by the %$\ggg$-invariant property of the form $(\cdot,\ccot)$ we have
%$$
%}
\end{observation}

\subsection{Compatibility with pairing $\langle\cdot,\cdot\rangle$}\label{sec: compatib}

Denote $V:=\bbr\otimes_\bbz \frakr$ and $V_0:=\bbr\otimes_\bbz \Phi$.
%
%
\iffalse
Then  $V_\bbz:= \bbz\frakr$ is the $\bbz$-form of $V$, which is also a $\bbz$-lattice in $\ttt^*$.
Now we already have a nondegenerate bilinear form $(\cdot,\cdot)_\ttt$ on $\ttt$. Furthermore, we choose the $\bbz$-form $\ttt_\bbz$ such that it is  the dual of $V_\bbz$, which means $V_\bbz=\ttt_\bbz^*$ because of the assumption that $\frakr$ spans $X(T)$.
Take a $\bbz$-basis $\{H_i\}$ in $\ttt_\bbz$. Then the dual basis $\{H_i^*\}$ in $\ttt^*$ spans a $\bbz$-basis of the corresponding $\bbz$-form $\ttt^*_\bbz$.
We can establish a nondegenerate bilinear form $(\cdot,\cdot)_{\ttt^*}$ on $\ttt^*$ via defining $(H_i^*,H_j^*)=(H_i,H_j)$. Through this process, we can establish a nondegenerate
symmetric bilinear form on $(\cdot,\cdot)_V$ on $V$.
 \fi
 %
 %
Then the following lemma shows that we can make a compatibility
normalization, % if necessary,
such that the bilinear form on
$(\cdot,\cdot)_{V_0}$ on $V_0$
and the good pairing $\langle\cdot,\cdot\rangle$
are compatible, that is,
they satisfy  BQR(3) of Definition \ref{defn: qrp}.

 \begin{lemma}\label{lem: compatibility normaliztion} Keep the notations and assumptions as above.
There is a bilinear form on $(\cdot,\cdot)_{V}$ on $V$ such that the related bilinear form $(\cdot,\cdot)$ on $\frakr$ is compatible with $\langle\cdot,\cdot\rangle$.
%, i.e. $\langle\cdot,\cdot\rangle){\Phi\times \Phi^\vee}$
%this is to say, both  satisfy
% BQR(3) of Definition \ref{defn: qrp}.
  \end{lemma}

 \begin{proof}
 For this,  we will  redefine a good pairing $\langle\cdot,\cdot\rangle$ on $X(T)\times X^\vee(T)$ if necessary such that
 $\langle \alpha,\delta^\vee\rangle=2(\delta,\delta)^{-1}(\alpha,\delta)$ for $\alpha,\delta\in \Phi$. The arguments are standard, the same as the one for  complex semisimple Lie algebras (see for example, \cite[\S8.3]{Hum1}), which is outlined below.

(Step 1) Note that $(\cdot,\cdot)_\ttt$ is  nondegenerate. For $\theta\in\frakr$, there is a unique $H_\theta$ such that $\theta(H)=(H_\theta,H)$  for all  $H\in\ttt$.
{We claim  that $[\ggg_\theta,\ggg_{-\theta}]=\bk H_\theta$. Actually, for any $X\in\ggg_{\theta}, Y\in \ggg_{-\theta}$, and $H\in \ttt$ we have
\begin{align*}
(H,[X,Y])&=([H,X], Y)  \cr
&=\theta(H)(X,Y)\cr
&=(H_\theta,H)(X,Y)\cr
&=((X,Y)H, H_\theta)\cr
&=(H, (X,Y)H_\theta).
\end{align*}
Hence $H$ is orthogonal to $[X,Y]-(X,Y)H_\theta$. By the nondegenerate property of $(\cdot,\cdot)_\ttt$ again,  we have
}
that for any $X\in \ggg_{\theta}$, $Y\in \ggg_{-\theta}$
\begin{align}\label{eq: H theta}
[X,Y]=(X,Y)H_\theta.
\end{align}
This confirms the claim. In particular,  $\alpha(H_\alpha)=(H_\alpha,H_\alpha)\ne 0$ for all $\alpha\in \Phi$.

(Step 2) For a generator $X_{\alpha}$ of $\ggg_{\alpha}$ for $\alpha\in \Phi$,  by the theory of complex semisimple Lie algebras there exists $X_{-\alpha}\in \ggg_{-\alpha}$ such that
$\{X_\alpha, h_\alpha, X_{-\alpha}\}$ forms an $\frak{sl}_2$-triple for $h_\alpha:=[X_\alpha, X_{-\alpha}]$, i.e. the Lie algebra generated by $X_{\pm\alpha}$ and $h_\alpha$ is Lie-algebra isomorphic to $\frak{sl}_2$ by  sending $X_\alpha, h_\alpha, X_{-\alpha}$ to $$  \left( \begin{array}{cc}
0 & 1\cr
0 & 0
\end{array}\right), \; \left( \begin{array}{cc}
1 & 0\cr
0 & -1
\end{array}\right), \; \left( \begin{array}{cc}
0 & 0\cr
1 & 0
\end{array}\right) \;
  $$
respectively. Furthermore, this $h_\alpha$ coincides with $2(H_\alpha,H_\alpha)^{-1} H_\alpha$, and $h_{-\alpha}=-h_\alpha$.

(Step 3) We can define for any $\beta,\theta\in \frakr$,
 \begin{align}\label{eq: sym form on roots-2}
 (\beta,\theta):=(H_\beta,H_\theta)
  \end{align}
 which gives rise to a symmetric bilinear form on $V$.

(Step 4) By the representation theory of $\frak{sl}_2$, $\delta(h_\alpha)\in \bbz$ for all $\alpha, \delta\in \Phi$, which is equal to  $2(\alpha,\alpha)^{-1}(\delta,\alpha)$.  So we can define
$\langle \delta,\alpha^{\vee}\rangle:=2(\alpha,\alpha)^{-1}(\delta,\alpha)\in\bbz$.
\end{proof}

 %By normalization if necessary, we can adjust the bilinear form such %that $(\cdot,\cdot)_{V_0}$ is compatible with the pairing %$\langle\cdot,\cdot\rangle$ as in
 %BQR(3) of Definition \ref{defn: qrp}. Furthermore, by normalization we %can make it sure that for $\gamma\in \Gamma$ and $\beta\in\frakr$,  %$\langle \beta,\gamma^\vee\rangle=2(\gamma,\gamma)^{-1}(\beta,\gamma)$ %if $(\gamma,\gamma)\ne 0$, and %$\langle\beta,\gamma^\vee\rangle=(\beta,\gamma)$ if otherwise.

\subsection{Two root arguments} Keep the notations and assumptions as above.
In this subsection we prove some auxiliary lemmas towards the properties
(S1) and (S2) in Def. \ref{defn: qrp} characterizing basic
quasi-reductive groups.

\begin{lemma}\label{lem: 4.2}
\begin{itemize}
\item[(1)] Let $\gamma, \nu\in \Gamma$ be linearly independent. Then
 $(\gamma, \nu) \ne0$ implies $\{\gamma\pm \nu\}\cap \frakr\ne
\emptyset$.
\item[(2)] For any $\gamma\in \Gamma$, if $(\gamma,\gamma)\ne 0$, then $2\gamma\in \Phi$.
%\commento{R. prefer to avoid extra notation $\frakr_0$}
  %=:\frakr_0$ where $\frakr_0$ denotes the even root system.
\end{itemize}
\end{lemma}

\begin{proof}
%The part (3) directly follows from the arguments in \S\ref{sec: %compatib}.
%
As to (1), consider the root subsystem generated by $\gamma,\nu$.
\iffalse

By (\ref{eq: gamma Lie bracket}), we have for $X_{\pm\gamma}\in \ggg_{\pm\gamma}$, there is a Cartan element
$H_\gamma=[X_\gamma, X_{-\gamma}]$.
 %with constant $c\in \bk^\times$ (see (\ref{eq: gamma Lie bracket})).
%Note that the characteristic of $\bk$ is zero.
So we have
%
\begin{align}\label{eq: comp two roots-1}
[X_\nu, [X_\gamma, X_{-\gamma}]]=[X_\nu, H_\gamma]=-\nu(H_\gamma)X_\nu.
\end{align}
%
On the other hand, by Jacobi identity we have
%
\begin{align}\label{eq: comp two roots}
[X_\nu, [X_\gamma, X_{-\gamma}]=[[X_\nu,X_\gamma], X_{-\gamma}]-[X_\gamma, [X_\nu, X_{-\gamma}]].
\end{align}
%
If $\{\gamma\pm \nu\}\cap \frakr= \emptyset$, equivalently,   $\{\nu\pm \gamma\}\cap \frakr= \emptyset$ (note that we already proved $\frakr=-\frakr$), then
the right hand side  of (\ref{eq: comp two roots}) is zero. Hence $\nu(H_\gamma)=0$, a contradiction.
\fi
%Then we can copy the arguments in the proof of Proposition
%\ref{prop: one-dimension} for (Case 2.2). So the part (1) is proved.

We need to show there exists one among $\gamma\pm\nu$ which belongs to $\frakr$.
By (\ref{eq: gamma Lie bracket}), we have that, for $X_{\pm\gamma}\in \ggg_{\pm\gamma}$, there is a Cartan element
$H_\gamma=[X_\gamma, X_{-\gamma}]$.
 So we have
\begin{align}\label{eq: comp two roots-1}
[X_\nu, [X_\gamma, X_{-\gamma}]]=[X_\nu, H_\gamma]=-\nu(H_\gamma)X_\nu.
\end{align}
On the other hand, by Jacobi identity we have
\begin{align}\label{eq: comp two roots}
[X_\nu, [X_\gamma, X_{-\gamma}]=[[X_\nu,X_\gamma], X_{-\gamma}]-[X_\gamma, [X_\nu, X_{-\gamma}]].
\end{align}
If $\{\gamma\pm \nu\}\cap \frakr= \emptyset$, equivalently,   $\{\nu\pm \gamma\}\cap \frakr= \emptyset$, then
the right hand side  of (\ref{eq: comp two roots}) is zero. Hence $\nu(H_\gamma)=0$,
contradicting the assumption $\gamma(H_\nu)\ne0$.

Now we prove (2), By the same reason as the above arguments for (1), if $2\gamma$ turns out to be a root, then it must be an even root. Now we show that $2\gamma\in \frakr$.
Consider
\begin{align}\label{eq: comp one roots-1}
[X_\gamma, [X_\gamma, X_{-\gamma}]]=[X_\gamma, H_\gamma]=-\gamma(H_\gamma) X_\gamma\ne 0.
\end{align}
On the other side,
\begin{align}\label{eq: comp one roots}
[X_\gamma, [X_\gamma, X_{-\gamma}]=[[X_\gamma,X_\gamma], X_{-\gamma}]-[X_\gamma, [X_\gamma, X_{-\gamma}]].
\end{align}
If $2\gamma\notin \frakr$, then $[X_\gamma,X_\gamma]=0$. Thus the right hand side of (\ref{eq: comp one roots}) becomes $-[X_\gamma, [X_\gamma, X_{-\gamma}]]$. Then
there would be $2[X_\gamma, [X_\gamma, X_{-\gamma}]]=0$, contradicting (\ref{eq: comp one roots-1}). Thus, (2) is proved.
\end{proof}

\begin{corollary}\label{coro: nonzero inn prod}
For $\gamma\in\Gamma$ with $(\gamma,\gamma)\ne0$,
$$\langle \beta,\gamma^\vee\rangle=2(\gamma,\gamma)^{-1}(\beta,\gamma)
\in\bbz$$ for $\beta\in\frakr$. And  $\textsf{r}_\gamma(\beta):=\beta-\langle \beta,\gamma^\vee\rangle \gamma\in\frakr$ for $\beta\in\frakr$.
\end{corollary}
\begin{proof} By Lemma \ref{lem: 4.2}(2), $\bar\gamma:=2\gamma\in \Phi$. Then we can make use of $\bar\gamma$ and the corresponding reflection $\textsf{r}_{\bar\gamma}\in \scrw$ {{because $\bar\gamma$ already becomes a root of the reductive group $G$}}. Note that $\ggg_\bo$ is a rational $G$-module, so the dimension of the weight space $\ggg_\gamma$ (=the odd root space) of $\ggg_\bo$ is preserved under $\scrw$-action, i.e. $\dim\ggg_{\gamma}=\dim\ggg_{\sigma(\gamma)}$ for any $\sigma\in\scrw$. {{Hence $\sigma(\gamma)\in\frakr$. Especially, $\textsf{r}_{\bar\gamma}(\beta)\in\frakr$ for $\beta\in \frakr$.}}  From this, the desired statement follows.
\end{proof}

\begin{remark} By Corollary \ref{coro: nonzero inn prod}, we can further require in Lemma \ref{lem: compatibility normaliztion} that for  $\gamma\in \Gamma$, we define
$$\langle\theta,\gamma^\vee\rangle=2(\gamma,\gamma)^{-1}(\theta,\gamma)$$
as long as $(\gamma,\gamma)\ne 0$.
 Clearly,
$\langle \beta_1, \beta_2^\vee\rangle \ne0$ if and only if  $\langle \beta_2, \beta_1^\vee\rangle \ne0$ for any $\beta_i\in \frakr$, $i=1,2$.
\end{remark}

%

%\subsubsection{}
\begin{lemma}\label{lem: 4.9} Let $\gamma\in \Gamma(D)$ with $(\gamma, \gamma)=0$. The following statements hold.
\begin{itemize}
\item[(1)] If $\gamma\pm \nu\in \frakr$ for $\nu\in \Gamma(D)$,
then $\gamma\pm \nu\in \Phi$. But $2\gamma\notin \Phi$.
\item[(2)] If $(\gamma, \beta)\ne 0$ for $\beta\in \frakr$, then $\{\gamma\pm\beta\}\cap \frakr\ne \emptyset$.
\item[(3)] Under the assumption in (2), {additionally $\beta\in\Gamma(D)$ with $(\beta,\beta)=0$,}  then
$\{\beta\pm k\gamma\mid k\in\bbz\}\cap \frakr\ne \emptyset$ implies
$k\in \{0, \pm1\}$.
     \end{itemize}
\end{lemma}

\begin{proof}
 %At first, according to the arguments in \S\ref{sec: compatib},  we have that for any $\beta_i\in \frakr$, $i=1,2$, the situation when $\langle \beta_1,\beta_2^\vee\rangle=0$ is equivalent to that of $(\beta_1,\beta_2)=0$.
 We proceed with verification by steps. %The part (1) is clear.
  The first part of (1) is because $\gamma\pm\nu\in \frakr$ implies that $\ggg_{\gamma\pm \nu}\ne 0$, and $[\ggg_\gamma, \ggg_{\pm\nu}]\subset \ggg_{\gamma\pm\nu}$. Hence the statement follows from the $\bbz_2$-gradation structure of $\ggg$.  As to the second part, we prove it by contradiction. Suppose $2\gamma\in\Phi$, then $(2\gamma,2\gamma)\ne0$. By the bilinear property of $(\cdot,\cdot)$ (see (Step 3) in the proof of Lemma \ref{lem: compatibility normaliztion}),  $(\gamma,\gamma)\ne0$, a contradiction.

Now we prove (2). It follows from Lemma \ref{lem: 4.2}(1) when $\beta\in\Gamma$. Suppose $\beta\in \Phi$, then by the weight chain property of $G$-rational module $\ggg_\bo$, $\textsf{r}_\beta(\gamma)=\gamma-\langle \gamma, \beta^\vee\rangle\beta$ is still a weight, and $\langle \gamma, \beta^\vee\rangle\ne0$ implies $\{\gamma\pm \beta\}\cap \Gamma\ne \emptyset$. The  part (2) is proved.

As to (3), it can be proved by the same arguments as in the proof of \cite[Lemma 1.10]{Ser1}.
\end{proof}

When $\gamma\in \Gamma$ with $(\gamma,\gamma)=0$, by Lemma \ref{lem: 4.9}(1),
we can define a map:
%\commento{R. shortened}
%$$\textsf{r}_\gamma: \frakr\rightarrow \frakr$$
%such that
% $\beta\in \frakr$

\begin{align}\label{eq: odd refl}
\textsf{r}_\gamma: \frakr\rightarrow \frakr, \qquad
\sfr_\gamma(\beta)=\begin{cases} \beta\pm \gamma &\text{ if } (\gamma,\beta)\ne 0;\cr
   \beta &\text{ if } (\gamma,\beta)= 0.
  \end{cases}
\end{align}

\iffalse
We immediately have the following consequence.
\begin{lemma}\label{lem: invertible}
The above  $\textsf{r}_\gamma$ is invertible.
\end{lemma}

\begin{proof}
For this, we only need to show that $\textsf{r}_\gamma^2=\id$. Arbitrarily take $\beta\in \frakr$.
If $(\gamma,\beta)=0$, then by definition $\textsf{r}_\gamma^2(\beta)=\beta$. Now suppose $(\gamma,\beta)\ne 0$. Then
$\textsf{r}_\gamma(\beta)\in \{\beta\pm \gamma\}$. Next $(\gamma,\beta\pm\gamma)=(\gamma,\beta)\ne0$. Hence
$\textsf{r}^2_\gamma(\beta)\in \{(\beta\pm\gamma)\mp \gamma\}=\{\beta\}$
 by Lemma \ref{lem: 4.9}(3).
\end{proof}
\fi

Furthermore, we have the following result, refining the
statement of Lemma \ref{lem: 4.9}(2).

\begin{lemma}\label{lem: only one from pair}
Let $\gamma\in \Gamma(D)$ with $(\gamma, \gamma)=0$. If $(\gamma, \beta)\ne 0$ for
$\beta\in \frakr$ and  {the odd reflection defined in (\ref{eq: odd refl}) is invertible,}  then
there is only one from the pair $\gamma\pm\beta$  which  belongs to $\frakr$.
\end{lemma}

\begin{proof} With Lemma \ref{lem: 4.9}(3), the statement can be proved by the same arguments as in the proof of  \cite[Lemma 1.11]{Ser1}. For the reader's convenience, we give a sketch.
  %proof.

The proof is by contradiction. Suppose both $\gamma\pm\beta$ are roots. If $\beta$ is an odd root, then $\gamma\pm\beta\in \Phi$. Then we have $(\gamma\pm\beta,\gamma\pm\beta)\ne0$ (see
(Step 1) in the proof of Lemma \ref{lem: compatibility normaliztion}). By Lemma \ref{lem: 4.9}(3), $\textsf{r}_\gamma(\beta+\gamma)=\beta$, and $\textsf{r}_\gamma(\beta-\gamma)=\beta$. This contradicts the invertible assumption  of $\textsf{r}_\gamma$.

Now suppose $\beta\in\Phi$. Then $\beta\pm\gamma\in\Gamma$. We first conclude that $\beta\pm\gamma$ must be non-isotropic root, i.e. $(\beta\pm\gamma,\beta\pm\gamma)\ne0$. Otherwise,  $(\beta\pm\gamma)\mp 2\gamma\in \frakr$ makes a contradiction with Lemma \ref{lem: 4.9}(3).
Thus, $(\beta\pm\gamma, \beta\pm\gamma)\ne 0$. By Corollary \ref{coro: nonzero inn prod}
with $\beta+\gamma$ and $\beta-\gamma$ respectively, we have $\langle \gamma, (\beta+\gamma)^\vee\rangle\in\bbz$, and $\langle \gamma, (\beta-\gamma)^\vee\rangle\in\bbz$.
%By the same arguments as in
%(\ref{eq: double integers result zero}),
By Corollary \ref{coro: nonzero inn prod} we have
 \begin{align*}
 &\langle \gamma, (\beta+\gamma)^\vee\rangle={{2(\gamma,\beta)}\over {(\beta,\beta)+2(\gamma,\beta)}}
 ={{\langle\gamma,\beta^\vee\rangle}\over{1+\langle \gamma,\beta^\vee\rangle}}\in\bbz, \cr
 &\langle \gamma, (\beta-\gamma)^\vee\rangle={{2(\gamma,\beta)}\over {(\beta,\beta)-2(\gamma,\beta)}}
 ={{\langle\gamma,\beta^\vee\rangle}\over{1-\langle \gamma,\beta^\vee\rangle}}\in\bbz
  \end{align*}
 which implies that $\langle \gamma, \beta^\vee\rangle=0$. This contradicts $(\gamma,\beta)\ne0$.
\end{proof}

\subsection{Main result}
 {
\begin{definition}
We say that a quasi-reductive supergroup is of \textit{monodromy type} if
$\bbG/\bbG_\ev \cong \bba^{(0|N)}$ where
$$
N=\#\Gamma(D).
$$
In other words in the decomposition (\ref{D-dec}) all $D_\gamma$ have dimension $1$.
\end{definition}

We are ready to state our main result.}

\begin{theorem}\label{thm: 4.1} Let $\bbG$ be a
connected quasi-reductive supergroup with %the purely-even part
$G:=\bbG_\ev$
and Lie superalgebra $\ggg=\Lie(\bbG)$ {{ satisfy the invertible property of all odd reflections defined in (\ref{eq: odd refl})}}.   Let $T$ be a maximal torus of  $G$ and
$\frakr=\Phi(G,T)\cup \Gamma(D)$ the root system of $\bbG$.
Suppose $\frakr \subset X(T)\backslash\{0\}$, $\frakr$  spans $X(T)$,
and $\ggg$ admits a nondegenerate even super-symmetric
$\ggg$-invariant bilinear form.  Then $\bbG$ is a basic quasi-reductive supergroup of monodromy type.
\end{theorem}

{{
\begin{remark}
In Kac's classification of finite-dimensional simple Lie
{superalgebras} over complex numbers (see \cite{Kac}), the Lie superalgebras with nondegenerate even super-symmetric bilinear form  coincide with the contragredient ones.   They are often called basic classical Lie superalgebras (see \cite[Remark 1.15]{CW}).    In the present paper, our notion is compatible with these references, { but $\widetilde{SL}(2|2)$  is excluded in our classification (Thm. \ref{main2})}.
\end{remark}
}}

In the following, we first show that under the hypothesis
of Theorem \ref{thm: 4.1}, all root spaces in the above
decomposition (\ref{eq: multiplicity ques}) are one dimensional.

\begin{proposition}\label{prop: one-dimension}
Let the notations and assumptions be as above.
All root spaces in the decomposition
$$
\ggg=\ttt\oplus \sum_{\alpha\in \Phi}\ggg_\alpha\oplus \sum_{\gamma\in \Gamma}\ggg_\gamma
$$
are one-dimensional.
\end{proposition}

\begin{proof}
%(3) Claim 1: $\dim\ggg_\gamma=1$ for $\gamma\in \Gamma$ with $\gamma(H_\gamma)=(H_\gamma,H_\gamma)\ne 0$. We prove this by reduction to absurdity.
Note that $G=\bbG_\ev$ is a connected reductive algebraic group with root system $\Phi$. Hence all root spaces $\ggg_\alpha$ for $\alpha\in \Phi$ are one-dimensional. So it suffices to prove the statement for $\ggg_\gamma$ with $\gamma\in \Gamma$.
By contradiction:
suppose $\dim\ggg_\gamma>1$. We then proceed with dividing the arguments into different cases.

(Case 1) Suppose $\gamma$ is a non-isotropic odd root, i.e. $(\gamma,\gamma)\ne 0$. { In this case,
 it is readily known that $[X_\gamma,X_\gamma]\ne0$. Actually, by (\ref{eq: gamma Lie bracket}) we have
 \begin{align}\label{eq: linear indp-1}
 [X_{-\gamma}, [X_\gamma,X_\gamma]]=2[H_\gamma,X_\gamma]=2\gamma(H_\gamma)X_\gamma\ne0
 \end{align}
 which implies $[X_\gamma,X_\gamma]$ is a nonzero vector of $\ggg_{2\gamma}$.  Hence $2\gamma\in \Phi$.  According to the previous arguments, $\ggg_{2\gamma}$ is one-dimensional. So $\ggg_{2\gamma}=\bk[X_\gamma,X_\gamma]$.

 }

 Under the assumption $\dim \ggg_\gamma>1$, it is readily shown that  there is another $X'_\gamma\in\ggg_\gamma$ linearly-independent of $X_\gamma$ such that
$(X'_\gamma, X_{-\gamma})=0$. {Actually, if $(X'_\gamma, X_{-\gamma})$ is not zero, say $c(\ne0)\in\bk$, then we take $X''={1\over c}X'_\gamma-X_{\gamma}$ which is still in $\ggg_{\gamma}$, linear-independently of  $X_{\gamma}$, and satisfies $(X''_\gamma, X_{-\gamma})=0$. This $X''_{\gamma}$ is desired. }

By (\ref{eq: H theta})
%the same arguments as for (\ref{eq: gamma Lie bracket}),
it is deduced that
\begin{align}\label{eq: pair of gamma gamma'}
[X'_\gamma,X_{-\gamma}]=0.
\end{align}
 {
Consider
$$X_1:=\ad(X_\gamma)(X'_{\gamma}).$$
Then $X_1\in \ggg_{2\gamma}$.
 By (\ref{eq: gamma Lie bracket})  and (\ref{eq: pair of gamma gamma'}), we have
 \begin{align}\label{eq: linear indp-2}
 [X_{-\gamma},X_1]=[H_\gamma,X_\gamma']=\gamma(H_\gamma)X'_\gamma\ne0.
 \end{align}
 Keep it in mind that $X_1$ and $[X_\gamma,X_\gamma]$ are two nonzero root vectors of $\ggg_{2\gamma}$, thereby linearly dependent.  Comparing (\ref{eq: linear indp-1}) and (\ref{eq: linear indp-2}), we can see that both left sides of both equations are linearly dependent. However $X'_\gamma$ and $X_\gamma$ are linearly independent, a contradiction.

 }

\iffalse
Generally, we have for $i\geq 1$
\begin{align*}
[X_{-\gamma}, X_i]&=[[X_{-\gamma},X_\gamma], X_{i-1}]-[X_\gamma, [X_{-\gamma}, X_{i-1}]]\cr
&=i \gamma(H_\gamma)X_{i-1}-[X_\gamma,[X_{-\gamma},X_{i-1}].
\end{align*}
Taking (\ref{eq: i=1}) into account,  we can inductively write  $[X_{-\gamma},X_{i}]=a_i(\gamma(H_\gamma))X_{i-1}$ with $a_i\in\bbz$ for $i\in \bbz_{>0}$. }
Consequently, we can prove
\begin{align}\label{eq: induction on i}
[X_{-\gamma}, X_i]=-((i-1)+a_{i-1})\gamma(H_\gamma)X_{i-1}
\end{align}
with $a_1=a_2=a_3=-1$, $\ldots, a_i=a_{i-1}-(i-1)\ne 0$ for $i\geq4$.
Note that $X_0=X'_\gamma\ne 0$ and already assume $\gamma(H_\gamma)\ne0$. So all $X_i$ are nonzero for $i\in\bbz_{\geq0}$. {\red {However, $X_2\in \ggg_{3\gamma}=0$, which is a contradiction.
\fi

(Case 2) Suppose $\gamma$ is isotropic, i.e. $(\gamma,\gamma)=0$.
 By the non-degeneracy of $(\cdot,\cdot)_\ttt$ again,  $\gamma(\ttt)\ne 0$. Note that $\frakr$ spans $X(T)$. So $\ttt$ can be spanned by $\{H_\theta\mid \theta\in \frakr\}$. So there exists $\theta\in \frakr$ such that $\gamma(H_\theta)\ne 0$ which means
 \begin{align}\label{eq: bilinear form of gamma theta nonzero}
 (\gamma,\theta)\ne0.
 \end{align}
 Thanks to Lemma \ref{lem: only one from pair}, there is only one from $\{\theta\pm \gamma\}$ belonging to $\frakr$. For example, $\theta-\gamma\in\frakr$, and then $\gamma-\theta\in\frakr$.
Consequently,
\begin{align}\label{eq: theta+gamma nonroot}
\pm(\theta+\gamma)\notin\frakr.
\end{align}

  We first assert that there exists a nonzero root vector $X''_\gamma\in\ggg_\gamma$ such that
  \begin{align}\label{eq: -theta gamma 0}
  [X_{-\theta}, X''_\gamma]=0
   \end{align}
   for a given root vector $X_{-\theta}\in \ggg_{-\theta}$. Actually, if $[X_{-\theta}, X_\gamma]=0$, then we choose $X''_\gamma=X_\gamma$. If not, then $0\ne X'_0:=[X_{-\theta},X_\gamma]\in \ggg_{\gamma-\theta}$.
   %Hence $\gamma-\theta\in \frakr$.
   As $\dim\ggg_\gamma>1$ by the beginning assumption, there exists  root vector $X_1\in \ggg_\gamma$ linearly independent of $X_\gamma$. If $X_1':=[X_{-\theta},{X_1}]\ne 0$, then $X_1'$ is a nonzero vector of $\ggg_{\gamma-\theta}$. If $X_0'$ and $X_1'$ are linearly dependent, then we easily choose a nonzero vector $X''_\gamma\in \ggg_{\gamma}$ such that $[X_{-\theta},X''_\gamma]=0$. If not, we can choose  $X_2\in\ggg_\gamma$ which is a nonzero proper linearly-span of $X_\gamma$ and $X_1$, satisfying $[X_{-\theta}, X_2]=0$. Then we take $X_\gamma''=X_2$.

As mentioned in (\ref{eq: gamma Lie bracket}), there exist $X_{-\gamma}''\in \ggg_{-\gamma}$ and $X_\theta\in\ggg_\theta$ such $[X_\gamma'',X_{-\gamma}'']=H_\gamma$ and
$[X_\theta,X_{-\theta}]=H_\theta$.  Hence we have the following
\begin{align*}
(H_\theta,H_\gamma)&=([X_\theta,X_{-\theta}], [X_\gamma'',X_{-\gamma}'']\cr
&=(X_\theta, [X_{-\theta},[X_\gamma'',X_{-\gamma}'']])\cr
&{\overset{(\ref{eq: -theta gamma 0})}{=}} (-1)^{\sfp(\theta)\sfp(\gamma)}(X_\theta,[X_\gamma'',[X_{-\theta},X_{-\gamma}'']])\cr
&=([X_\theta,X_\gamma''],[X_{-\theta},X_{-\gamma}'']),
\end{align*}
where $\sfp(\theta)$ and $\sfp(\gamma)$ denote the parity of the roots $\theta$ and $\gamma$, respectively.
Note that both $\pm(\theta+\gamma)$ are not roots. Hence the final part of the above equations is equal to zero.  This contradicts the assumption $(H_\gamma,H_\theta)=(\gamma,\theta)\ne0$ in (\ref{eq: bilinear form of gamma theta nonzero}).

By the same arguments as above we can deal with the case  when $\pm (\theta+\gamma)\in \frakr$ and $\pm(\theta-\gamma)\notin\frakr$.

Summing up, we accomplish the proof.
\end{proof}

\iffalse
\begin{remark} The above proposition (and its proof) shows that there impossibly exists any nondegenerate  super-symmetric invariant bilinear form on $\tilde \ggg$ constructed in Remark \ref{remark:  monodromy}.
\end{remark}
\fi

\subsection{Proof of Theorem \ref{thm: 4.1}}  Due to Proposition \ref{prop: one-dimension},
%and Theorem \ref{thm: 3.11},
it suffices to show that $\frakr=\Phi\cup \Gamma$ satisfies the axioms in Definition \ref{defn: qrp}.
By the arguments in \S\ref{sec: compatib}, we only need to verify BQR(4) Serganova's condition. More precisely, we only need to verify the second parts (S1) and (S2) of BQR(4), respectively.
%By the arguments in \S\ref{sec: compatib}, one knows that  %$\langle\gamma,\beta^\vee\rangle=0$ is equivalent to $(\gamma,\beta)=0$.
%
When $\gamma\in \Gamma$ with $(\gamma,\gamma)=0$, by Lemma \ref{lem: 4.9}(1) we can define a mapping
$$\textsf{r}_\gamma: \frakr\rightarrow \frakr$$
such that
 $\beta\in \frakr$
 \begin{align*}
 \sfr_\gamma(\beta)=\begin{cases} \beta\pm \gamma &\text{ if } (\gamma,\beta)\ne 0;\cr
   \nu &\text{ if } (\gamma,\beta)= 0.
  \end{cases}
\end{align*}
{{By assumption, this $\textsf{r}_\gamma$ is invertible.}}  Thus, the axiom BQR(4) (S1) is satisfied.

As to BQR(4)(S2), it follows from Corollary \ref{coro: nonzero inn prod}.
 According to the analysis in the beginning of the proof, we already proved the theorem.

{

\begin{remark}\label{rem: excluding sl(2_2)-case}
In Theorem  \ref{thm: 4.1}, we assume that all  odd reflections defined in (\ref{eq: odd refl}) are invertible. This assumption does not meet the case $A(1|1)$. We demonstrate this point with an example.

  Let $\mathbb{G}$ be a supergroup of type $A(1|1)$. Recall that $\mathbb{G}$ admits the standard fundamental root system $\{\epsilon_1-\epsilon_2, \epsilon_2-\delta_1, \delta_1-\delta_2\}$ and the root system $\frak{R}=\Phi\cup \Gamma$ with
 $\Phi=\{\pm \epsilon_1\mp \epsilon_2, \pm \delta_1\mp \delta_2\}$ and $\Gamma=\{\pm\epsilon_i\mp \delta_j\mid i,j=1,2\}$ subject to the relation $\epsilon_1+\epsilon_2-(\delta_1+\delta_2)=0$. There is indeed a nondegenerate bilinear form on $\ggg=\text{Lie}(\mathbb{G})$ which gives rise to the bilinear form on the root system $\frak{R}$ via defining $(\epsilon_i,\epsilon_j)=\delta_{ij}$,  $(\delta_i,\delta_j)=-\delta_{ij}$ where $i,j=1,2$, and $\delta_{ij}=1$ if $i=j$ and $0$ otherwise.
 We take $\beta=\epsilon_2-\epsilon_1$, and $\gamma=\delta_1-\epsilon_2$. Then $\beta\in \mathfrak{R}$ and $\gamma\in \Gamma$ satisfy $(\gamma,\gamma)=0$ and $(\beta,\gamma)\ne0$. Then we have  the following facts
 \begin{align*}
 \beta &\in \mathfrak{R};\cr
 \beta+\gamma =\delta_1-\epsilon_1  &\in \mathfrak{R}; \text{ and}\cr
 \beta+2\gamma=\delta_1-\delta_2   &\in \mathfrak{R}.
 \end{align*}
If defining the odd reflection $\sfr_\gamma$ as in (\ref{eq: odd refl}), then we have $\sfr_\gamma(\beta)=\beta+\gamma$ and $\sfr_\gamma(\beta+2\gamma)=\beta+\gamma$. Hence $\sfr_\gamma$ is not invertible.

   So the supergroup $\mathbb{G}$ of type $A(1|1)$ is surely not a basic quasi-reductive supergroup because it violates  the axiom (S1) of BQR(4) in Definition \ref{defn: qrp}.
   Actually,  a quasi-reductive supergroup  does not contain any component subgroup isomorphic to the one of type $A(1|1)$. This explains the reason why type $A(1|1)$ is {{not contained }} in  the  list of the classification theorem of indecomposable basic quasi-reductive supergroups (Theorem \ref{thm: 2.13}).
\end{remark}
}

\end{document}